\documentclass[12pt,leqno]{article}
\usepackage{amsfonts}
\usepackage{amsmath, amsthm, amsfonts, amssymb, color}
\usepackage{mathrsfs}
\usepackage{color}
\theoremstyle{definition}

\newcommand{\scr}[1]{\mathscr #1}
\definecolor{wco}{rgb}{0.5,0.2,0.3}

\numberwithin{equation}{section} \theoremstyle{remark}

\newcommand{\ua}{\uparrow}

\title{{\bf    Exponential Ergodicity for Time-Periodic  McKean-Vlasov SDEs }\footnote{ Supported in
 part by  the National Key R\&D Program of China (No. 2022YFA1006000, 2020YFA0712900), NNSFC (11921001) and 
  ERC Advanced Grant.} }
\author{
{\bf Panpan Ren$^{a)}$,  Karl-Theodor Sturm$^{b)}$,  Feng-Yu Wang$^{c)}$  }\\
\footnotesize{$^{a)}$ Mathematics department, Hong Kong City University}\\
\footnotesize{$^{b)}$ Mathematics department, Boon University,  Germany}\\
\footnotesize{$^{c)}$ Center for Applied Mathematics, Tianjin University, Tianjin 300072, China}\\
\footnotesize{ rppzoe@gmail.com, sturm@iam.uni-bonn.de,  wangfy@tju.edu.cn}}
\begin{document}
\allowdisplaybreaks
\def\R{\mathbb R}  \def\ff{\frac} \def\ss{\sqrt} \def\B{\mathbf
B}
\def\N{\mathbb N} \def\kk{\kappa} \def\m{{\bf m}}
\def\ee{\varepsilon}\def\ddd{D^*}
\def\dd{\delta} \def\DD{\Delta} \def\vv{\varepsilon} \def\rr{\rho}
\def\<{\langle} \def\>{\rangle}
  \def\nn{\nabla} \def\pp{\partial} \def\E{\mathbb E}
\def\d{\text{\rm{d}}} \def\bb{\beta} \def\aa{\alpha} \def\D{\scr D}
  \def\si{\sigma} \def\ess{\text{\rm{ess}}}\def\s{{\bf s}}
\def\beg{\begin} \def\beq{\begin{equation}}  \def\F{\scr F}
\def\Ric{\mathcal Ric} \def\Hess{\text{\rm{Hess}}}
\def\e{\text{\rm{e}}} \def\ua{\underline a} \def\OO{\Omega}  \def\oo{\omega}
 \def\tt{\tilde}\def\[{\lfloor} \def\]{\rfloor}
\def\cut{\text{\rm{cut}}} \def\P{\mathbb P} \def\ifn{I_n(f^{\bigotimes n})}
\def\C{\scr C}      \def\aaa{\mathbf{r}}     \def\r{r}
\def\gap{\text{\rm{gap}}} \def\prr{\pi_{{\bf m},\varrho}}  \def\r{\mathbf r}
\def\Z{\mathbb Z} \def\vrr{\varrho} \def\ll{\lambda}
\def\L{\scr L}\def\Tt{\tt} \def\TT{\tt}\def\II{\mathbb I}
\def\i{{\rm in}}\def\Sect{{\rm Sect}}  \def\H{\mathbb H}
\def\M{\mathbb M}\def\Q{\mathbb Q} \def\texto{\text{o}} \def\LL{\Lambda}
\def\Rank{{\rm Rank}} \def\B{\scr B} \def\i{{\rm i}} \def\HR{\hat{\R}^d}
\def\to{\rightarrow}\def\l{\ell}\def\iint{\int}\def\gg{\gamma}
\def\EE{\scr E} \def\W{\mathbb W}
\def\A{\scr A} \def\Lip{{\rm Lip}}\def\S{\mathbb S}
\def\BB{\scr B}\def\Ent{{\rm Ent}} \def\i{{\rm i}}\def\itparallel{{\it\parallel}}
\def\g{{\mathbf g}}\def\Sect{{\mathcal Sec}}\def\T{\mathcal T}\def\BB{{\bf B}}
\def\f{\mathbf f} \def\g{\mathbf g}\def\BL{{\bf L}}  \def\BG{{\mathbb G}}
\def\Bd{{D^E}} \def\BdP{D^E_\phi} \def\Bdd{{\bf \dd}} \def\Bs{{\bf s}} \def\GA{\scr A}
\def\Bg{{\bf g}}  \def\Bdd{\psi_B} \def\supp{{\rm supp}}\def\div{{\rm div}}
\def\ddiv{{\rm div}}\def\osc{{\bf osc}}\def\1{{\bf 1}}\def\BD{\mathbb D}\def\GG{\Gamma}
\def\H{{\bf H}} \def\n{{\bf n}}
\maketitle

\begin{abstract} As extensions to the corresponding results derived for time homogeneous Mckean-Vlasov SDEs, the exponential ergodicity is proved for time-periodic  distribution dependent SDEs in three different situations:
  \beg{enumerate}
 \item[1)]  in the quadratic Wasserstein distance and relative entropy
for the dissipative case;
\item[2)]   in the Wasserstein distance induced by a cost function for the partially dissipative case; and
\item[3)]  in the weighted Wasserstein distance induced by a cost function and a Lyapunov function  for the fully non-dissipative case.
\end{enumerate}
The main results are illustrated by time inhomogeneous  granular media equations, and are extended to reflecting McKean-Vlasov SDEs in a convex domain. \end{abstract} \noindent

 \vskip 2cm 

 \section{Introduction}

Recently, by using the log-Harnack and Talagrand inequalities,  the exponential ergodicity in relative entropy is proved in \cite{20RW}  for a class of McKean-Vlasov SDEs, which include as typical examples the granular porous media equations
investigated in \cite{CMV,GLW}. Next, by using coupling methods, the exponential ergodicity in different probability metrics have been derived in \cite{W21a} for partially dissipative and non-dissipative models. Moreover, these types of exponential ergodicity  have been investigated in \cite{W21b} for reflecting McKean-Vlasov SDEs. In this paper, we  extend these results to time-periodic (reflecting) McKean-Vlasov SDEs.

Let $D\subset \R^d$ be a convex domain. When $D\ne \R^d$, it has a non-empty boundary $\pp D$. In this case, for any $x\in\pp D$ and $r>0$, let
$$\scr N_{x,r}:=\big\{\n\in \R^d: |\n|=1, B(x-r\n, r)\cap D=\emptyset\big\},$$ where $B(x,r):=\{y\in \R^d: |x-y|<r\}.$ We have
$$\scr N_x:=\cup_{r>0} \scr N_{x,r}\ne\emptyset, \ \ x\in \pp D, r>0.$$ We call $\scr N_x$ the set of inward unit normal vectors of $\pp D$ at point $x$. Since $D$ is convex, $\scr N_x\ne\emptyset$ for $x\in\pp D$ and
\beq\label{CVX} \<x-y,\n(x)\>\le 0,\ \ y\in \bar D, x\in \pp D, \n(x)\in \scr N_x.\end{equation}
Let $\scr P(\bar D)$ be the space of all probability measures on the closure $\bar D$ of $D$,   equipped with the weak topology.
 Consider the following reflecting McKean-Vlasov SDE on $\bar D\subset \R^d$:
\beq\label{E1} \d X_t= b_t(X_t,\L_{X_t}) \d t+ \si_t(X_t,\L_{X_t})\d W_t + \n(X_t) \d l_t,\ \ t\ge 0,\end{equation}
 where   $W_t$ is an $m$-dimensional Brownian motion on a complete filtration probability space $(\OO,\{\F_t\}_{t\ge 0},\P)$, $\L_{X_t}$ is the distribution of $X_t$,
 $\n(x)\in\scr N_x$ for $x\in \pp D$, $l_t$ is an adapted increasing process which increases only when $X_t\in \pp D$, and
$$b: [0,\infty)\times \R^d\times \scr P(\bar D) \to \R^d,\ \  \si: [0,\infty)\times \R^d\times \scr P(\bar D)\to \R^d\otimes\R^m$$ are measurable.
When $D=\R^d$ we simply denote $\scr P=\scr P(\bar D)$. In this case,  we have $\pp D=\emptyset$ so that $l_t=0$ and \eqref{E1} reduces to
\beq\label{E01} \d X_t= b_t(X_t,\L_{X_t}) \d t+ \si_t(X_t,\L_{X_t})\d W_t,\ \ t\ge 0.\end{equation}

The SDE \eqref{E1} or \eqref{E01} is called well-posed for distributions in a subspace $\hat{\scr P}\subset \scr P(\bar D)$, if  for any $s\ge 0$ and any $\F_s$-measurable variable $X_s$ with $\L_{X_s}\in \hat{\scr P}$, \eqref{E1} has a unique solution $(X_t)_{t\ge s}$ with $\L_{X_\cdot}\in C([s,\infty);\hat{\scr P})$, the space of continuous maps from $[s,\infty)$ to $\hat{\scr P}$ under the weak topology. In this case, we denote
$P_{s,t}^*\mu=\L_{X_t}$ for  the solution with $\L_{X_s}=\mu\in \hat{\scr P}$.

In this paper, we investigate the exponential ergodicity of \eqref{E1} and \eqref{E01} with   $t_0$-periodic coefficients for some $t_0>0:$
$$(b_{t+t_0},\si_{t+t_0})=(b_t, \si_t),\ \ t\ge 0,$$  such that the corresponding results derived in \cite{20RW, W21a, W21b} are extended to time inhomogeneous models. By the $t_0$-periodicity and the well-posedness for distributions in $\hat {\scr P},$  we have
\beq\label{PP1} P_{s, t}^*\mu= P_{s+nt_0, t+nt_0}^*\mu,\ \ t\ge s\ge 0, n\in \mathbb N, \mu\in \hat{\scr P}.\end{equation} In this case, a probability measure $\bar \mu_0\in \hat{\scr P}$ is called an invariant probability measure, if
$P_{0, t_0}^* \bar\mu_0= \bar\mu_0.$  Combining this with \eqref{PP1}, we see   that the measures
$$\bar\mu_s:= P_{0,s}^*\bar\mu_0,\ \ s\in [0,t_0]$$
satisfy
\beq\label{PP2} P_{s+m t_0, s+(m+n) t_0}^*\bar \mu_s= \bar\mu_s,\ \ n,m\in \mathbb Z_+, s\in [0,t_0].\end{equation}
Let $\W: \hat{\scr P} \times \hat{\scr P} \to [0,\infty)$ with $\W(\mu,\nu)=0$ if and only if $\mu=\nu$. We call \eqref{E1} exponential ergodic
in $\W,$ if there exist constants $c,\ll>0$ such that $P_t^*:=P_{0,t}^*$ satisfies
\beq\label{ECC} \W(P_{s,s+nt_0}^*\mu,\bar\mu_s)\le c \e^{-\ll n} \W(\mu,\bar\mu_s),\ \ n\in\mathbb N, \mu\in \hat{\scr P}.\end{equation}
By \eqref{PP1}, this is equivalent to
$$ \W(P_{s+mt_0, s+(m+n)t_0}^*\mu,\bar\mu_s)\le c \e^{-\ll n} \W(\mu,\bar\mu_s),\ \ n, m\in\mathbb Z_+, \mu\in \hat{\scr P}, s\in [0,t_0].$$
So, we will only consider \eqref{ECC}.

The remainder of the paper is organized as follows. In Sections 2-4, we study the exponential ergodicity  for \eqref{E01} without reflection, where   Section 2   considers dissipative models for  $\W$ being  the quadratic Wasserstein distance $\W_2$ or the relative entropy $\H$,  Section 3 concerns with partially dissipative models with $\W=\W_\psi$ induced by a cost function $\psi$,    and   Section 4 deals with fully non-dissipative models for $\W=\W_{\psi,V}$   induced by a cost function $\psi$ and a Lyapunov function $V$.
Finally,  these results are extended in Section 5 to  the reflecting SDE \eqref{E1} on a convex domain   $D$.

\section{Exponential ergodicity in  relative entropy and $\W_2$}

Corresponding to   \cite{CMV,GLW,20RW} where the exponential ergodicity in entropy is investigated in the time homogeneous case,
we consider the exponential ergodicity   in relative entropy for \eqref{E01}. Recall that  the relative entropy  for   probability measures $\mu_1,\mu_2\in \scr P$ is given by
$$\H(\mu_1|\mu_2):= \beg{cases} \mu_2(\rr\log \rr), &\text{if} \ \rr:=\ff{\d\mu_1}{\d\mu_2}\ {\rm exists},\\
\infty, &\text{otherwise.} \end{cases}$$
For  the symmetric diffusion process generated by $L:=\DD+\nn V$ on $\R^d$ with $\bar\mu(\d x):= \e^{V(x)}\d x\in \scr P$,
the exponential ergodicity in $\H$ with rate $\ll>0$ is equivalent to the   log-Sobolev inequality
$$\bar\mu(f^2\log f^2) \le \ff 2 \ll \bar\mu(|\nn f|^2) ,\ \   f\in C_b^1(\R^d), \bar\mu(f^2)=1,$$
where  $\mu(f):=\int f\d\mu$ for a measure $\mu$ and $f\in L^1(\mu)$.
  According to the concentration property of the log-Sobolev inequality (see \cite{AMS}), there exists $\vv>0$ such that $\bar\mu(\e^{\vv |\cdot|^2})<\infty$,
  so that by Young's inequality, $\H(\mu|\bar\mu)<\infty$ implies
  $$\mu(|\cdot|^2)\le \vv^{-1} \big\{\H(\mu|\bar\mu)+ \log \bar\mu(\e^{\vv|\cdot|^2})\big\}<\infty.$$
  Therefore,    to investigate the exponential convergence in entropy, it is natural to consider distributions in the Wasserstein space
$$\scr P_2:= \big\{\mu\in \scr P: \mu(|\cdot|^2)<\infty\big\},$$
which is  a Polish space  under the quadratic Wasserstein distance
$$\W_2(\mu_1,\mu_2):= \inf_{\pi\in\C(\mu_1,\mu_2)} \bigg(\int_{\R^d\times\R^d} |x-y|^2\pi(\d x,\d y)\bigg)^{\ff 1 2},\ \ \mu_1,\mu_2\in \scr P_2,$$
where $\C(\mu_1,\mu_2)$ is the set of all couplings of $\mu_1$ and $\mu_2$. 

\subsection{Assumptions}

 Let $\dd_x$ be the Dirac measure at $x\in\R^d$.
We assume
\beg{enumerate} \item[$(H_1)$]     $|b_t(0,\dd_0)|+\|\si_t(0,\dd_0)\|$ is locally integrable in $t\ge 0$,  and  there exist  $K_1,K_2,K_3\in L_{loc}^1( [0,\infty);\R)$
 such that
 \beg{align*}   &\|\si_t(x,\mu)-\si_t(y,\nu)\|^2\le K_3(t)\big(|x-y|^2+\W_2(\mu,\nu)^2\big),\\
 &2\<b_t(x,\mu)-b_t(y,\nu),x-y\>+\|\si_t(x,\mu)-\si_t(y,\nu)\|_{HS}^2\\
  & \le K_1(t) |x-y|^2 + K_2(t) \W_2(\mu,\nu)^2,\   \ t\ge 0, x,y\in \R^d, \mu,\nu\in
\scr P_2.\end{align*}
\end{enumerate} According to \cite[Theorem 3.3]{20HRW} (see also \cite{W18}), under this condition the SDE \eqref{E01} is well-posed for distributions in $\scr P_2$, and
\beq\label{EX0} \W_2(P_{s,t}^*\mu,P_{s,t}^*\nu)^2\le \e^{\int_s^t (K_1(r)+K_2(r))\d r} \W_2(\mu,\nu)^2,\ \ t\ge s\ge 0, \mu,\nu\in \scr P_2.\end{equation}
To deduce from \eqref{EX0} the exponential ergodicity in entropy, we need the following condition.

\beg{enumerate} \item[$(H_2)$]  $\si_t(x,\mu)=\si_t(x)$ does not depend on $\mu$ and is invertible, and
 there exist  increasing  positive measurable functions  $\ll,\kk_1,\kk_2$ such that
\beg{align*} &2\<b_t(x,\mu)-b_t(y,\nu),x-y\>^++\|\si_t(x)-\si_t(y)\|_{HS}^2
   \le \kk_1(t) |x-y|^2 + \kk_2(t) |x-y| \W_2(\mu,\nu),\\
&\|\si_t(x)^{-1}\| \le \ll(t),\ \ t\ge 0, x,y\in\R^d, \mu,\nu\in \scr P_2.\end{align*}   \end{enumerate}
Obviously, $(H_2)$ implies    $(H_1)$ for $K_1(t)=\kk_1(t)+\bb_t$ and $K_2(t)=\ff{\kk_2(t)^2}{4\bb_t}$ for $\bb_t>0$, but in applications we may take better choices of $(K_1,K_2)$   than that implied by $(H_2)$. For any $t\ge s\ge 0$, let
$$\ll(s,t):= \sup_{r\in [s,t]} \ll(r),\ \ \kk_i(s,t):= \sup_{r\in [s,t]} \kk_i(r),\ \ i=1,2.$$

We intend to establish the following type of estimate
\begin{equation}\label{EST}
\H(P^*_{s,t}\mu|P^*_{s,t}\nu)\leq \phi(s, t)\H(\mu|\nu),\\\ t>s, \mu\in \scr P_2
\end{equation}
for a reasonable class of measures $\nu\in \scr P_2.$ In the time homogeneous situation, one takes $\nu$ as the invariant probability measures so that $\P^*_{s,t}\nu=\nu$ for all $t\geq s.$

As explained above, to derive \eqref{EST}, we  need to establish the log-Sobolev inequality for $P_{s,t}^*\nu$. To this end, we apply the Bakry-Emery curvature for the associated  time-distribution dependent generator of \eqref{E01}:
$$L_{t,\mu}:= \ff 1 2 {\rm tr}\big\{\si_t\si_t^*\nn^2\big\}+ b_t(\cdot,\mu)\cdot\nn,\ \ t\ge 0, \mu\in \scr P_2.$$
According to \cite{CM}, we introduce
\beq\label{GGG} \beg{split} &\GG_{t}^{1}(f,g):= \ff 1 2 \<\si_t\si_t^*\nn f,\nn g\>,\ \ f,g\in C^1(\R^d),\\
&\GG_{t,\mu}^2(f,f):= \ff 1 2 L_{t,\mu} \GG_t^1(f,f)- \GG_t^1(f, L_{t,\mu}f)+\ff 1 2 \pp_t \GG_t^1(f,f),\ \ f\in C^3(\R^d).\end{split}\end{equation}
To make $\GG_{t,\mu}^2$ meaningful and also for late use, we assume

\beg{enumerate} \item[$(H_3)$]   $A_t:=\|\si_t\|_\infty \in L_{loc}^2([0,\infty)),$ at least one of the following two conditions holds:
\item[$(1)$] $\si_t$ is constant for each $t\ge 0$;
\item[$(2)$] $\si_t(x)$ is   $C^1$ in $t$ and $C^2$ in $x$, $b_t(x)$ is $C^1$ in $x$, and there exists a function $\gg\in L_{loc}^1([0,\infty);\R)$ such that
$$\GG_{t,\mu}^2(f,f) \ge \gg_t\, \GG_t^1(f,f),\ \ t\ge 0, f\in C^3(\R^d), \mu\in \scr P_2.$$\end{enumerate}

Finally, for any constant $c>0$, we write $\nu\in T_c$ if $\nu\in \scr P$ satisfying the Talagrand inequality
\beq\label{TNN} \W_2(\mu, \nu)^2\le  c \H(\mu|\nu).\end{equation}
According to \cite{BGL}, this inequality is implied by the log-Sobolev inequality
\beq\label{LNN} \nu(f^2\log f^2) \le c \nu(|\nn f|^2),\ \ f\in C^1_b(\R^d), \nu(f^2)=1,\end{equation}
for which we denote $\nu\in Log_c$.

 \subsection{Main results}

\beg{thm}\label{T1.1} Assume $(H_1)$ and  that $\eqref{E01}$ is $t_0$-periodic for some $t_0>0$ with
\beq\label{LMM} \ll:=- \int_0^{t_0} \{K_1(r)+K_2(r)\}\d r>0.\end{equation}
\beg{enumerate} \item[$(1)$]  $\eqref{E01}$ has a unique invariant probability measure $\bar \mu_0$ such that
\beq\label{EXW} \W_2(P_{nt_0}^* \mu, \bar\mu_0)^2\le \e^{-n \ll} \W_2(\mu,\bar\mu_0)^2,\ \ \mu\in\scr P_2, n\in\mathbb N.\end{equation}
\item[$(2)$] If $(H_2)$, and one of  $(H_3)(1)$ or $(H_3)(2)$   with $\int_0^{t_0} \gg_s\d s>0$ hold, then there exists a constant $c>0$ such that for any $n\in\mathbb N$ and $\mu\in \scr P_2,$
\beq\label{EXM} \max\big\{\H(P_{nt_0}^*\mu|\bar\mu_0), \W_2(P_{nt_0}^* \mu, \bar\mu_0)^2\big\}
  \le c \e^{-\ll n} \min\big\{\H(\mu|\bar\mu_0),  \W_2(\mu,\bar\mu_0)^2\big\}. \end{equation}
\end{enumerate}
\end{thm}

To illustrate this result, we consider the time-dependent version of granular media equations studied in \cite{ CMV,GLW,20RW}. Let $V\in C^{0, 2}([0,\infty)\times\R^d)$ and $W\in C^{0,2}([0,\infty)\times \R^{2d})$ such that
\beq\label{VW}  \int_{\R^d} \e^{-V_t(x)} \d x+\int_{\R^d\times\R^d} \e^{-V_t(x)-V_t(y)-\ll W_t(x,y)}\d x\d y <\infty, \ \ \ll>0,t\ge 0.\end{equation} Consider the following PDE on $\D_2$, the space of all probability density functions on $\R^d$ such that the corresponding probability measure is in $\scr P_2$:
\beq\label{PDEW} \pp \rr_t= {\rm div}\big\{\nn \rr_t-\rr_t \nn(V_t+ W_t\circledast \rr_t) \big\},\end{equation}
where for a probability measure $\mu$ or a probability density function $\rr$
$$W_t\circledast \mu:=\int_{\R^d} W_t(\cdot,y)\mu(\d y),\ \ W_t\circledast \rr:= \int_{\R^d} W_t(\cdot,y)\rr(y)\d y.$$
We will use $\nn^{(1)}$ and $\nn^{(2)}$ to denote the gradient operators in the first and second components on the product space $\R^d\times \R^d,$ so that
$$\|\nn^{(1)}\nn^{(2)} W_t(x,y)\| := \sup_{u,v\in\R^d, |u|,|v|\le 1}  |\nn_u^{(1)}\nn_v^{(2)} W_t(x,y)|,\ \ t\ge 0, x,y\in\R^d,$$ where $\nn_u$ stands for the directional derivative along $u$.
We let
$$\|\nn^{(1)}\nn^{(2)} W_t\|_\infty:=\sup_{x,y\in\R^d} \|\nn^{(1)}\nn^{(2)} W_t(x,y)\|.$$
For any probability density $\rr$ on $\R^d$ and any $s\ge 0$, let $P_{s,t}^*\rr$ be the solution of \eqref{PDEW} for $t\ge s$ and $\rr_s=\rr$.
If $(V_t,W_t)$ is $t_0$-periodic, $\bar \rr_0\in\D_2$ is called an invariant solution of \eqref{VW} if $P_{0,t_0}^*\bar \rr_0=\bar\rr_0$. In this case, let
$$\bar\rr_s:= P_{0,s}^* \bar\rr_0,\ \ s\in (0,t_0).$$
Moreover, for any two probability density functions $\rr_1,\rr_2$,
$$\H(\rr_1|\rr_2):= \H(\rr_1(x)\d x|\rr_2(x)\d x).$$
 Let $I_d$ be the $d\times d$ identity matrix.

   \beg{thm}\label{T1.2}  Let $(V_t,W_t)$ be $t_0$-periodic for some $t_0>0$, and  there exists   $\gg\in L_{loc}([0,t_0]; \R)$   such that $\ll:=\int_0^{t_0}\gg_t\d t>0$ and
 \beq\label{CV} \Hess_{V_t+W_t(\cdot,z)} \ge (\gg_t + \|\nn^{(1)}\nn^{(2)}W_t\|_\infty)I_d,\ \ t\in [0,t_0], z\in \R^d.\end{equation}
Then $\eqref{PDEW}$ has a unique invariant solution $\bar\rr_0$ such that
 \beq\label{EXM'}\beg{split} &\max\big\{\W_2(\rr_{nt_0}(x)\d x, \bar\rr_0(x)\d x)^2, \H(P_{nt_0}^*\rr|\bar\rr_0)\big\}\\
 &\le c \e^{-\ll  n} \min\big\{\W_2(\rr_{0}(x)\d x, \bar\rr_0(x)\d x)^2,\H(\rr|\bar\rr_0)\big\},\ \ n\in\mathbb N,\rr\in\D_2.\end{split}\end{equation}
 \end{thm}

 \subsection{Proofs  }
We first prove the following lemma which also applies to the non-periodic case.
\beg{lem}\label{LN1} Assume $(H_1)$, $(H_2).$
For any $t\ge s\ge 0$,    let
\beq\label{PHI} \phi_{s,t}:= \ll(s,t)^2 \Big(\ff{\kk_1(s,t)}{1-\e^{-\kk_1(s,t)}} +\ff{(t-s)\kk_2(s,t)^2} 2 \e^{2(t-s)\kk_1(t)+2\kk_2(t)}\Big).\end{equation}
\beg{enumerate} \item[$(1)$] For any $\vv>0, c>0$ and $\nu\in T_c$,
\beq\label{EXX1}\beg{split}  &\H(P_{s,t}^*\mu|P_{s,t}^*\nu)\le \phi_{t-\vv,t} \e^{\int_s^{t-\vv} (K_1+K_2)(r)\d r} \W_2(\mu,\nu)^2\\
&\le c \phi_{t-\vv,t}   \e^{\int_s^{t -\vv} (K_1+K_2)(r)\d r}\H(\mu,\nu),\ \ t\ge s+\vv, s\ge 0, \mu\in \scr P_2.\end{split}\end{equation}
\item[$(2)$]  If $(H_3)(2)$ holds and $\nu\in Log_c$ for some constant $c>0$, then
\beq\label{EXX2} \beg{split} &\H(P_{r,t}^*\mu|P_{s,t}^*\nu)\le \phi_t \e^{\int_r^{t-\vv} (K_1+K_2)(r)\d r} \W_2(\mu,P_{s,r}^*\nu)^2\\
&\le c(s,r) \phi_{t-\vv,t} \e^{\int_r^{t-\vv} (K_1+K_2)(r)\d r}\H(\mu|P_{s,r}^*\nu),\ \ t\ge r+\vv, r\ge s\ge 0, \mu\in \scr P_2\end{split} \end{equation} holds for
$$c(s,r):=  cA_r^2 \ll(s,r)^2   \e^{-2\int_s^r \gg_\theta\d\theta}+ 4A_r^2\int_s^r\e^{-2\int_\tau^r\gg_\theta\d \theta}\d \tau,\ \ r\ge s\ge 0.$$
\end{enumerate} \end{lem}

\beg{proof}  (1) By   \cite[Theorem 4.1]{W18} or \cite[Theorem 4.1]{20HRW},  assumption $(H_2)$ implies
 $$\H(P_{s,t}^* \mu|P_{s,t}^*\nu)\le \phi_{s,t} \W_2(\mu,\nu)^2,\ \ t\ge s\ge 0, \mu,\nu\in \scr P_2.$$
 So, for $t\ge s+\vv$ we obtain
 \beq\label{LH}\H(P_{s,t}^* \mu|P_{s,t}^*\nu)= \H(P_{t-\vv,t}^* P_{s,t-\vv}^*\mu| P_{t-\vv,t}^* P_{s,t-\vv}^*\nu) \le \phi_{t-\vv,t }\W_2(P_{s,t-\vv}^*\mu, P_{s,t-\vv}^*\nu)^2.\end{equation}
 Next, by \cite[Theorem 3.1]{W18},  assumption $(H_1)$ implies
$$\W_2(P_{s,t-\vv}^*\mu, P_{s,t-\vv}^*\nu)^2\le \e^{2\int_s^{t-\vv} (K_1(r)+K_2(r))\d r} \W_2(\mu,\nu)^2.$$
 Combining this with \eqref{LH} and applying \eqref{TNN} we prove \eqref{EXX1}.

(2)  Noting that $P_{s,t}^*\nu =P_{r,t}^*(P_{s,r}^*\nu)$, to deduce \eqref{EXX2} from \eqref{EXX1} we need only to prove
 $P_{s,r}^*\nu\in T_{c(s,r)}$ which follows from  $P_{s,r}^*\nu\in Log_{c(s,r)}.$ To this end, we let $\nu_t:= P_{s,t}^*\nu$ and consider the decoupled (classical) SDE of \eqref{E01}:
 \beq\label{DCP} \d X_t^\nu= b_t(X_t^\nu, \nu_t) \d t + \si_t(X_t^\nu)\d W_t,\ \ t\ge s, X_0^\nu\in \R^d.\end{equation}
 For any $\mu\in \scr P$, let $(P_{s,t}^\nu)^*\mu=\L_{X_{s,t}^\nu}$ for $\L_{X_s^\nu}=\mu$. Then
 \beq\label{MKK} P_{s,r}^*\nu= (P_{s,t}^\nu)^* \nu, \ \ r\ge s.\end{equation}
 Now, for $\nu\in Log_c,$    $\|\si_s^{-1}\|\le \ll_s$ implies
 \beq\label{MM1} \nu(f\log f)\le \ff{c \ll_s^2}4 \nu\Big(\ff{|\si_s^* \nn f|^2}f\Big),\ \ 0<f\in C_b^1(\R^d), \nu(f)=1. \end{equation}
 According to \cite[Theorem 4.1]{CM} for the time inhomogeneous Markov semigroup associated with \eqref{DCP}, we remark that in this result $\GG(f)$ is misprint from $\ff{\GG(f)}f$ (see Lemma \ref{LN2} below for $D=\R^d$),  $(H_3)$ and \eqref{DCP} yield  that $\nu_r:= P_{s,r}^*\nu= (P_{s,t}^\nu)^*\nu$ satisfies
\beq\label{MM2} \nu_r(f\log f)\le \ff{c(s,r) }{4A_r^2} \nu_r\Big(\ff{|\si_r^* \nn f|^2}f\Big),\ \ 0<f\in C_b^1(\R^d), \nu_r(f)=1.\end{equation}
 Sine $\|\si_r\|_\infty\le A_r$,  this  implies $P_{s,r}^*\nu\in T_{c(s,r)}$ as desired.

\end{proof}

 \beg{proof}[Proof of Theorem \ref{T1.1}]  By shifting a time $s\in [0,t_0),$ for simplicity, we assume $s=0.$\\
 (1) By \eqref{LMM} and the $t_0$-periodicity,   the uniqueness of $\bar\mu$ and \eqref{EXW} follows from \eqref{EX0}. So, it suffices to prove the existence of $\bar\mu_0$.

 Take
\beq\label{RPP} \mu_n:= P_{0, nt_0}^*\dd_0,\ \ n\in\mathbb N.\end{equation}
 We intend to prove that $\mu_n$ converges to some $\bar\mu_0\in \scr P_2$ as $n\to\infty$, so that by a standard argument   the semigroup property
 of $\bar P_n^*:= P_{0,n t_0}^*$:
 $$\bar P_{n+m}^*= \bar P_n^*\bar P_m^*,\ \ n,m\in \mathbb Z_+,$$
 we conclude that $\bar\mu_0$ is an invariant probability measure. To this end, it remains to show that $\{\mu_n\}_{n\ge 1}$ is a $\W_2$-Cauchy sequence, i.e.
 \beq\label{CAU} \lim_{n\to\infty}\sup_{k\ge 1} \W_2(\mu_n, \mu_{n+k})=0.\end{equation}
 By  \eqref{EX0},  \eqref{PP1} and \eqref{LMM},  we obtain
 \beq\label{W1} \W_2(\mu_n, \mu_{n+k})^2\le \e^{-\int_0^{nt_0} (K_1(r)+K_2(r))\d r} \W_2(\dd_0, P_{0,k t_0}^*\dd_0)^2=\e^{-\ll n} \E |X_{kt_0}|^2,\end{equation}
 where $X_t$ solves \eqref{E01} with $X_0=0.$ By taking $y=0, \nu=\dd_0$ in $(H_1)$, and noting that the periodicity and $(H_2)$ implies that $|b_\cdot(0,\dd_0)|+\|\si_t(0)\|$ is bounded
 and $\|\si_t(x)\| \le c_0(1+|x|)$ for some constant $c_0>0,$ we find   constants $c_1,c_2>0$ such that
 \beg{align*} &2 \<b_t(x,\mu), x\> +\|\si_t(x)\|_{HS}^2\\
 &= 2\<b_t(x,\mu)-b_t(0,\dd_0), x-0\>+ \|\si_t(x)-\si_t(0)\|_{HS}^2\\
 &\quad + 2\<b_t(0,\dd_0), x\>-\|\si_t(0)\|_{HS}^2+2\<\si_t(x),\si_t(0)\>_{HS} \\
 &\le K_1(t)|x|^2 +K_2(t)\mu(|\cdot|^2) + c_1(1+|x|)\\
 &\le c_2 +\Big(K_1(t)+\ff\ll{2 t_0} \Big) |x|^2 +K_2(t) \mu(|\cdot|^2),\ \ t\ge 0, x\in\R^d, \mu\in\scr P_2.\end{align*}
 So, by applying It\^o's formula to \eqref{E01} for $X_0=0$, we obtain
 $$\d |X_t|^2\le \Big\{ \Big(c_2+K_1(t)+\ff\ll{2 t_0} \Big) |X_t|^2 +K_2(t) \E|X_t|^2\Big\}\d t+ \d M_t$$
 for some martingale $M_t$. By Duhamel's formula,   this and $X_0=0$ implies
\beq\label{SSA} \E|X_t|^2\le c_2 \int_0^t \e^{\int_s^t (K_1(r)+K_2(r)+\ff{\ll}{2 t_0})\d r} \d s,\ \ t\ge 0.\end{equation}
 By \eqref{LMM}  and the $t_0$-periodicity,  we find a constant $C>0$ such that
 $$\int_{s}^{s+kt_0}  \Big(K_1(r)+K_2(r)+\ff{\ll}{2 t_0}\Big)\d r = -\ff{\ll k} 2<0,\ \ s\ge 0, k\in\mathbb Z_+.$$
 So, letting $\lfloor r\rfloor:=\sup\{n\in\mathbb Z_+: r\ge n\}$ for $r\ge 0$, \eqref{SSA} implies
 \beq\label{SOP} \beg{split} \sup_{t\ge 0} \E|X_t|^2 &\le \sup_{t\ge 0} c_2 \int_0^t \e^{-\ff{\ll \lfloor(t-s)/t_0\rfloor}2 +\int_0^{t_0}|K_1(r)+K_2(r)+\ff{\ll}{2 t_0}|\d r}\d s\\
 &\le C \int_0^t \e^{-\ff{(t-s)\ll }2}\d s\le \ff{2C}{\ll}<\infty.\end{split}\end{equation}
 Combining this with \eqref{W1}, we prove the desired \eqref{CAU}.

 (2)    By \eqref{EXX2} and \eqref{EXW}, it suffices to find a constant $c>0$ such that $\bar\mu_0\in Log_c.$

 a) When $(H_3)(1)$ holds, let $(\bar P_{s,t})_{t\ge s}$ be the semigroup associated with the SDE
\beq\label{SDEB} \d \bar X_t = b_t(\bar X_t,\bar\mu_0)\d t + \si_t\d W_t,\end{equation}
that is,  letting $(\bar X_{s,t}^x)_{t\ge s}$ be the solution starting from $x$ at time $s$,
$$\bar P_{s,t}f(x):= \E f(\bar X_{s,t}^x),\ \ f\in \B_b(\R^d), t\ge 0.$$
By $(H_1)$ which implies $K_2\ge 0$,    we have
\beq\label{POL}\beg{split} &2\<b_t(x,\bar\mu_0)-b_t(y,\bar\mu_0), x-y\>\le K_1(t)|x-y|^2,\ \ x,y\in\R^d, t\ge 0,\\
&\int_0^{t_0} K_1(s)\d s\le -\ll<0.\end{split}\end{equation}  Then  $P_t:=P_{0,t}$ satisfies
\beq\label{GRD} |\nn \bar P_{s,t} f|\le \e^{\ff 1 2 \int_s^tK_1(s)\d s }\bar P_{s,t}|\nn f|\le c_1 \e^{-\ff \ll 2 \lfloor(t-s)/t_0\rfloor} \bar P_t|\nn f|,\ \ t\ge s\ge 0\end{equation}
for some constant $c_1>0.$  So, for any $f\in C_b^1(\R^d)$,
\beg{align*} &\bar P_t(f^2\log f^2)- (\bar P_t f^2)\log (\bar P_t f^2)= \int_0^t \ff{\d}{\d s} \bar P_s \big\{(\bar P_{s,t}f^2)\log (\bar P_{s,t}^* f^2)\big\}\d s\\
&= \int_0^t \bar P_{s} \ff{|\si_s^*\nn \bar P_{s,t} f^2|^2}{\bar P_{s,t} f^2}\d s\le c_1^2\int_0^t  \|\si_s\|^2 \e^{-\ll\lfloor(t-s)/t_0\rfloor} \bar P_s\bar P_{s,t}|\nn f|^2\d s\\
&= (\bar P_t|\nn f|^2)  c_1^2\int_0^t  \|\si_s\|^2 \e^{-\ll \lfloor(t-s)/t_0\rfloor} \d s,\ \ t\ge 0.\end{align*}
By the $t_0$-periodicity and $\|\si_\cdot\|^2\in L^1([0,t_0])$, we obtain
\beg{align*} & \int_0^t  \|\si_s\|^2 \e^{-\ll \lfloor(t-s)/t_0\rfloor} \d s
\le \sum_{i=0}^{\lfloor t/t_0\rfloor}\int_{it_0}^{(i+1)t_0} \|\si_s\|^2\e^{-\ll(\lfloor t/t_0\rfloor-i-1)}\d s\\
&\le \bigg(\int_0^{t_0}\|\si_s\|^2\d s\bigg)\sum_{i=0}^\infty \e^{-(i-1)\ll}=:c<\infty,\ \ t\ge 0,\end{align*}
so that there exists a constant $c>0$ such that
$$\bar P_t(f^2\log f^2)- (\bar P_t f^2)\log (\bar P_t f^2)\le c \bar P_t|\nn f|^2,\  \ t\ge 0, f\in C_b^1(\R^d).$$
 Moreover, by  \eqref{POL},
\eqref{SDEB} is exponential ergodic with unique invariant probability measure $\bar \mu_0$ as it reduces to \eqref{E01} when
$\L_{\bar X_0}=\bar\mu_0$. By taking $t=nt_0$ and letting $n\to\infty$, we prove $\mu_n\in Log_c$ for all $n\ge 1.$

b) When $(H_3)(2)$ holds with $\gg:=\int_0^{t_0}\gg_s\d s>0$,   we apply Lemma \ref{LN1} for  $s=0$ and $\nu=\dd_0$. Then \eqref{MM1} holds for $c=0$, so that by the $t_0$-periodic  and $\gg>0$, we find a constant $c'>0$ such that
$$c(0,nt_0)= 4\dd_{t_0}^2\int_0^{nt_0} \e^{-2\int_\tau^r
\gg_\theta\d \theta}\d \tau\le c',\ \ n\in\mathbb N.$$ Moreover,
by \eqref{MM2}, $\mu_n:= P_{0,n t_0}^*\dd_0$ satisfies
 $$ \mu_n(f^2\log f^2)\le \ff{c(0,nt_0) }{\|\si_{t_0}\|_\infty^2} \mu_n( |\si_{t_0}^* \nn f|^2)\le c(0, nt_0) \mu_n(|\nn f|^2),\ \ 0<f\in C_b^1(\R^d), \mu_n(f^2)=1.$$
 Therefore,   $\mu_n\in Log_{c'}$ for all $n\ge 1$, which together with \eqref{EXW} implies $\bar\mu_0\in Log_{c'}$.
\end{proof}

 \beg{proof}[Proof of Theorem \ref{T1.2}]
  It is easy to see  that for any $s\ge 0$ and probability density function $\rr$, $P_{s,t}^*\rr$ is the   density function of $\L_{X_t}$ for
$X_t$ solving \eqref{E01} from time $s$ with $\L_{X_s}=\rr(x)\d x$ and
 \beq\label{SB} \si_t(x):=\ss{2}I_d,\ \ b_t(x,\mu):= -\nn\big\{V_t+W_t\circledast\mu\big\}(x),\ \ t\ge 0, x\in]\R^d, \mu\in \scr P_2.\end{equation}
 Then \eqref{CV} implies
  \beg{align*} & 2\<b_t(x,\mu)-b_t(y,\nu), x-y\>= 2\<b_t(x,\mu)-b_t(y,\nu), x-y\> + 2 \<b_t(y,\mu)-b_t(y,\nu), x-y\>\\
& = - 2 \int_{\R^d} \mu(\d z) \int_0^1 \big\<\Hess_{V_t+W_t(\cdot,z)}(x + r (y-x)) (x-y),  x-y\big\> \d r\\
&\quad  + 2 \big\<\nu(\nn^{(1)} W_t(y,\cdot))- \mu(\nn^{(1)} W_t(y,\cdot)),x-y\big\>\\
&\le -2\big(\gg_t+\|\nn^{(1)}\nn^{(2)} W_1\|_\infty)|x-y|^2 + 2\|\nn^{(1)}\nn^{(2)} W_t\|_\infty  |x-y| \W_1(\mu,\nu)\\
&\le -2 \big(\gg_t+\|\nn^{(1)}\nn^{(2)} W_1\|_\infty)|x-y|^2 + 2\|\nn^{(1)}\nn^{(2)} W_t\|_\infty  |x-y| \W_2(\mu,\nu).\end{align*}
 Thus, $(H_1)$ holds for
 \beq\label{KKT}  K_1(t)= - 2\gg_t-\|\nn^{(1)}\nn^{(2)} W_t\|_\infty,\ \ K_2(t)= \|\nn^{(1)}\nn^{(2)} W_t\|_\infty,\end{equation}
 and $(H_2)$ holds for
 \beq\label{kkt} \kk_1(t)= -2 \big(\gg_t+\|\nn^{(1)}\nn^{(2)} W_t\|_\infty),\ \ \kk_2(t)= 2  \|\nn^{(1)}\nn^{(2)} W_t\|_\infty.\end{equation}
 Moreover, since $\si$ is constant, $(H_3)(1)$ holds.
Therefore,    this result follows from    Theorem \ref{T1.1}.
 \end{proof}

 \section{Ergodicity for partially dissipative models}

For any $\psi\in \Psi=\{\psi\in C^2([0,\infty)): \psi(0)=0, \psi'>0, \|\psi'\|_\infty<\infty\},$   let
$$\scr P_{\psi}:=\big\{\mu\in \scr P: \mu(\psi(|\cdot|))<\infty\big\},$$
$$\W_\psi(\mu,\nu):= \inf_{\pi\in \C(\mu,\nu)} \int_{\R^d\times \R^d} \psi(|x-y|)\pi(\d x,\d y),\ \ \mu,\nu\in \scr P_\psi.$$
Then  $\scr P_{\psi}$ is complete under $\W_\psi$, i.e. a $\W_\psi$-Cauchy sequence in $\scr P_\psi$ converges with respect to $\W_\psi.$
Let $\|\cdot\|_{Lip}$ be the Lipschitz constant for functions on $\R^d$. We assume

\beg{enumerate}
\item[$(H_4)$](Ellipticity) $\si_t(x,\mu)=\si_t(x)$ does not depend on $\mu$,   and  there exist    $\aa\in L_{loc}^1([0,\infty); (0,\infty))$ and a measurable map
$$\hat \si: [0,\infty)\times \R^d\to \R^d\otimes\R^d$$ such that
$$\sup_{t\in [0,T]}\big\{\|\si_t\|_{Lip}+ \|\hat \si_t\|_{Lip}\big\}<\infty,\ \ T>0,$$
$$\si_t(x)\si_t(x)^*=\aa_t    I_d+ \hat \si_t(x)\hat \si_t(x)^*,\ \ t\ge 0, x\in\R^d.$$
\item[$(H_5)$] (Partial dissipativity) Let  $\psi\in \Psi$, $\gg\in C([0,\infty))$ with $\gg_t(r)\le K r$ for some constant $K>0$ and all $r\ge 0$,   such that
\beq\label{A2E}  2\aa_t \psi''(r) +(\gg_t\psi')(r)\le -\kk_t \psi(r),\ \ r\ge 0\end{equation} holds for some $\kk\in L_{loc}^1([0,\infty);\R).$
Moreover, $b$ is   bounded on bounded subsets of $[0,\infty)\times \R^d\times \scr P_\psi$, and   there exists  $\theta\in L_{loc}^1([0,\infty);(0,\infty))$  such that
\beq\label{A3E} \beg{split}  & \<b_t(x,\mu)-b_t(y,\nu), x-y\> +\ff 1 2 \|\hat \si_t(x)-\hat \si_t(y)\|_{HS}^2\\
&\quad \le |x-y|  \big\{\theta_t  \W_\psi (\mu,\nu)  + \gg_t(|x-y|)\big\},\ \ t\ge 0, x,y\in \R^d, \mu,\nu\in \scr P_\psi.   \end{split}\end{equation}\end{enumerate}

\beg{thm}\label{T3-1} Assume  $(H_4)$  and $(H_5),$  with $\psi''\le 0$ if $\hat \si_t$ is non-constant for some $t\ge 0$. Then   $\eqref{E01}$ is well-posed with distributions in
 $\scr  P_\psi$, and   $P_t^*$ satisfies
 \beq\label{EXP1'0}  \W_\psi(P_t^*\mu, P_t^*\nu)\le  \e^{-\int_0^t\{\kk_s-\theta_s\|\psi'\|_\infty\}\d s}  \W_\psi(\mu,\nu),\ \ t\ge 0, \mu,\nu \in \scr P_\psi.\end{equation}
 Consequently,
 if   $(b_t,\si_t)$ is $t_0$-periodic, $\psi'(t)\le C\psi'(s)$ for some constant $C>1$ and all $t\ge s\ge 0$, and
 $$\ll:= \int_0^{t_0}   \{\kk_s-\theta_s\|\psi'\|_\infty\}\d s>0,$$
 then $\eqref{E01}$ has a unique invariant probability measure
$\bar\mu_0\in \scr P_\psi$ such that
\beq\label{EXP2'0}  \W_\psi(P_{s,s+nt_0}^*\mu, \bar\mu_0)\le  \e^{-n\ll}  \W_\psi(\mu,\bar\mu_0),\ \ n\in \mathbb N, \mu \in \scr P_\psi.\end{equation}
  \end{thm}

\beg{proof} The well-posedness and \eqref{EXP1'0}  follow from \cite[Theorem 3.1]{W21a} by using coupling methods. So, it suffices to prove the existence of the invariant probability measure $\bar\mu_0$ when  $\ll>0$ and the coefficients are $t_0$-periodic.  Let $x_0\in \R^d$. It suffices to show that the sequence $\{P_{nt_0}^*\dd_{x_0}\}_{n\ge 1}$ is a $\W_\psi$-Cauchy sequence so that its limit is an invariant probability measure of $\eqref{E01}$. By \eqref{EXP1'0} we have
$$\W_\psi(P_{nt_0}^*\dd_{x_0}, P_{(n+m)t_0}^*\dd_{x_0})\le C \e^{-n\ll} \W_\psi(\dd_{x_0}, P_{m t_0}^*\dd_{x_0}),\ \ n,m\ge 1.$$
Since $\ll>0$, it suffices to prove
\beq\label{PRT} \sup_{m\ge 1} \W_\psi(\dd_{x_0}, P_{m t_0}^*\dd_{x_0})<\infty.\end{equation}
By $\psi'(t)\le C\psi'(s)$ for $t\ge s$, we have
$$\psi(s+t)-\psi(s)=\int_{s}^{s+t}\psi'(r)\d r\le C\int_0^t \psi'(r)\d r= C\psi(t),\ \ s,t\ge 0.$$
This implies
$$\psi\Big(\sum_{i=1}^n s_i\Big)\le C\sum_{i=1}^n \psi(s_i),\ \ s_i\ge 0, n\ge 1.$$
Consequently,  by \eqref{EXP1'0} and $\ll>0$, we obtain
$$\W_\psi(\dd_{x_0}, P_{nt_0}^*\dd_{x_0})\le C\sum_{i=0}^{n -1}\W_\psi(P_{it_0}^*\dd_{x_0}, P_{(i+1)t_0}^*\dd_{x_0})\le C \W_\psi(\dd_{x_0}, P_{t_0}^*\dd_{x_0})\sum_{i=0}^\infty \e^{-i\ll}<\infty.$$
Therefore, \eqref{PRT} holds.
 \end{proof}

To illustrate Theorem \ref{T3-1}, we present below an example associated with  time-inhomogeneous granular media equations.
Let    $\W_1=\W_\psi$ and $\scr P_1(\R^d)= \scr P_\psi(\R^d)$ for $\psi(r)=r$.

\paragraph{Example 3.1.} Let $\aa\in L^1([0,t_0]:   (0,\infty))$ and
$$V: [0,t_0]\times \R^d\to\R,\ \ W: [0,t_0]\times\R^d\times\R^d\to\R$$ be measurable
with $V_t\in C^2(\R^d), W_t\in C^2(\R^d\times \R^d),$ and for some  constants $R, \theta_1,\theta_2>0$,
\beq\label{HDD} \Hess_{V_t+W_t(\cdot,z)}\ge \big(\theta_2\aa_t 1_{\{|\cdot|>R/2\}}- \theta_1\aa_t 1_{\{|\cdot|\le R/2\}}\big)I_d,\ \
t\in [0,t_0], z\in \R^d.\end{equation}
Consider \eqref{E01} with $t_0$-periodic coefficients
\beq\label{HDD2}\si_t=\ss{\aa_t}I_d,\ \ b_t(x,\mu):= -\nn \{V_t+W_t \circledast\mu\}(x),\ \ (t,x,\mu)\in [0,t_0]\times\R^d\times\scr P_1.\end{equation}
Let  $\gg(r)= \theta_1(r\land R)- \theta_2(r-R)^+$ for $r\ge 0$, and
$$\psi(r):=\int_0^r \e^{-\int_0^s\gg(u)\d u}\d s\int_s^\infty t\e^{\int_0^t\gg(u)\d u}\d t,\ \ r\ge 0.$$
Then
$$c_1(\psi):= \inf_{r\ge 0}\psi'(r)>0,\ \ c_2(\psi):=\sup_{r\ge 0} \psi'(r)<\infty.$$
If
$$ \ll:= 2\int_0^{t_0} \Big(\ff{\aa_t}{c_2(\psi)} -\ff{\|\nn^{(1)}\nn^{(2)}W_t\|_\infty}{c_1(\psi)}\Big)\d t>0,$$
then \eqref{E01} has a unique invariant probability measure $\bar\mu_0\in \scr P_1$ such that
$$\W_\psi(P_{nt_0}^*\mu,\bar\mu_0)\le   \e^{-\ll n} \W_\psi(\mu,\bar\mu_0),\ \ n\in \mathbb N, \mu\in \scr P_1=\scr P_\psi.$$ Consequently,
$$\W_1(P_{nt_0}^*\mu,\bar\mu_0)\le \ff{c_2(\psi)}{c_1(\psi)} \e^{-\ll n} \W_1(\mu,\bar\mu_0),\ \ n\in \mathbb N, \mu\in \scr P_1.$$

\beg{proof} It is easy to see that  $\psi\in C^2([0,\infty)$ with $\psi'>0$ and
$$\lim_{r\to\infty}\psi'(r) =\lim_{r\to\infty} \ff{\int_r^\infty t \e^{\int_0^t\gg(u)\d u}\d t}{\e^{\int_0^r \gg(u)\d u}}= \lim_{r\to\infty}\ff{r}{-\gg(r)}=\ff 1 {\theta_2}.$$
So, $0<c_1(\psi)<c_2(\psi)<\infty$ and
$$c_1(\psi)\W_1\le \W_\psi\le c_2(\psi)\W_1.$$
By Theorem \ref{T3-1}, it suffices to verify \eqref{A2E} and \eqref{A3E} for
\beq\label{GMM} \gg_t(r):= 2\aa_t \gg(r),\  \ \ \theta_t:=\ff 2 {c_1(\psi)} \|\nn^{(1)}\nn^{(2)}W_t\|_\infty,\ \ \kk_t:= \ff{2\aa_t}{c_2(\psi)}.\end{equation}

Firstly, by the definitions of $\gg$ and $\psi$, $\gg_t:= 2\aa_t\gg$ in \eqref{GMM}   we have
$$2\aa_t\psi''(r)+ 2\aa_t(\psi'\gg)(r)= -2\aa_tr\le -\ff{2\aa_t}{c_2(\psi)} \psi(r),\ \ r\ge 0.$$
Then \eqref{A2E} holds for $\gg_t$ and $\kk_t$ in \eqref{GMM}.

Next, by \eqref{HDD} and \eqref{HDD2}, we have $\hat\si=0$, and as in the proof of Theorem \ref{T1.2},
\beg{align*}& 2\<b_t(x,\mu)-b_t(y,\nu), x-y\>\\
&= -2 \int_{\R^d} \mu(\d z) \int_0^1 \<\Hess_{V_t+W_t(\cdot,z)}(x+r(x-y)) (x-y), x-y\>\d r\\
 &\quad + 2 \<\nu(\nn^{(1)}W_t(y,\cdot))-\mu(\nn^{(1)}W_t(y,\cdot)), x-y>\\
&\le 2 |x-y|^2 \int_0^1 \big(\theta_1\aa_t 1_{\{|x+r(y-x)|\le R/2\}} -\theta_2\aa_t1_{\{|x+r(y-x)|\ge R/2\}}\big)\d r\\
 &\quad + 2 \|\nn^{(1)}\nn^{(2)}W_t\|_\infty \W_1(\mu,\nu)|x-y|\\
&\le 2 \aa_t |x-y|\gg(|x-y|)+ \ff 2 {c_1(\psi)} \|\nn^{(1)}\nn^{(2)}W_t\|_\infty |x-y|\W_\psi(\mu,\nu)\end{align*}
holds for any $t\in [0,t_0], \ x,y\in\R^d$ and $ \mu,\nu\in \scr P_{\psi}=\scr P_1.$
Hence, \eqref{A3E} holds for $\gg_t$ and $\theta_t$ in \eqref{GMM}.
\end{proof}

\section{Ergodicity for non-dissipative models}

We consider the fully non-dissipative case such that \cite[Theorem 2.1]{W21a} is extended to the periodic  setting.
For any $t\ge 0$ and $\mu\in \scr P$, consider the second-order differential operator
\beq\label{LM'} L_{t,\mu}:= \ff 1 2 {\rm tr}\{\si_t\si_t^*\nn^2\}+ b_t(\cdot,\mu)\cdot\nn.  \end{equation}
 For any positive measurable function $V$ on $\R^d$, let
 $$\scr P_V:=\{\mu\in \scr P: \mu(V)<\infty\}.$$

\beg{enumerate} \item[$(H_7)$] (Lyapunov Condition) There exist    $0\le V\in C^2(\R^d)$ with $\lim_{|x|\to\infty} V(x) =\infty$ and
$K_0,K_1\in L^1_{loc}([0,\infty);\R)$ such that
\beq\label{H11}  \sup_{t\ge 0; x\in\R^d}
\ff{ |\si(t,x)\nn V(x)|}{1+V(x)} <\infty, \end{equation}
\beq\label{H120} L_{t,\mu} V\le K_0(t)-K_1(t)V,\ \ t\ge 0, \mu\in \scr P_V. \end{equation}\end{enumerate}

Since $\lim_{|x|\to\infty} V(x)=\infty,$ \eqref{H120} controls the long distance  behaviour of the associated stochastic system.
To ensure the exponential ergodicity, we also need conditions in short distance.
For any $l>0$, consider the    class
$$\Psi_l:= \big\{\psi\in C^2([0,l]; [0,\infty)):\ \psi(0)=\psi'(l)=0, \psi'|_{[0,l)}>0 \big\}.$$
For each $\psi\in \Psi_l,$ we extend it to the half line by setting $\psi(r)=\psi(r\land l)$, so that $ \psi'$ is  non-negative and  Lipschitz continuous with compact support and
\beq\label{HP} c_\psi:= \sup_{r>0} \ff{r\psi'(r)}{\psi(r)} <\infty. \end{equation}
For any  constant $\bb>0$, define the quasi-distance on $\scr P_V(\R^d)$:
$$\W_{\psi,\bb V} (\mu,\nu):= \inf_{\pi\in \scr C(\mu,\nu)} \int_{\R^d\times\R^d} \psi (|x-y|) \big(1+\bb V(x)+\bb V(y)\big)\pi(\d x,\d y),\ \ \mu,\nu\in \scr P_V.$$
To prove the exponential convergence of $P_t^*$ under $\W_{\psi,\bb V}$,  the dependence on distribution for the drift will be characterized by
\beq\label{HW}\beg{split}
 \hat \W_{\psi,\bb V}(\mu,\nu)&:=\inf_{\pi\in \C(\mu,\nu)} \ff{\int_{\R^d\times\R^d} \psi(|x-y|) (1+\bb V(x)+\bb V(y))  \pi(\d x, \d y)}{\int_{\R^d\times\R^d}  \psi'(|x-y|) (1+\bb V(x)+\bb V(y)) \pi(\d x, \d y)}\\
 &\ge \ff{\W_{\psi,\bb V}(\mu,\nu)}{\|\psi'\|_\infty (1+\bb \mu(V)+\bb \nu(V))},\ \ \mu,\nu\in \scr P_V.\end{split}\end{equation}

\beg{enumerate} \item[$(H_8)$] (Local monotonicity)  $b$ is bounded on bounded set in $[0,\infty)\times \R^d\times\scr P_V$. Moreover,    there exist      $l>0, \psi\in \Psi_l$ and $u_l,K,\theta\in L_{loc}^1([0,\infty); [0,\infty))$  such that
$$2\aa_t \psi''(r)+ K(t) \psi'(r) \le  - u_l (t)\psi(r), \ \ r\in [0,l],t\ge 0,$$
 \beg{align*}  & \<b_t(x,\mu)-b_t(y,\nu), x-y\> +\ff 1 2 \|\hat \si_t(x)-\hat\si_t(y)\|_{HS}^2\\
&\le K_t |x-y|^2 +\theta_t |x-y|\hat\W_{\psi,\bb V}(\mu,\nu),\ \ x,y\in \R^d, \mu,\nu \in \scr P_V,t\ge 0.   \end{align*}
\end{enumerate}

By $(H_7)$,  for any $l>0$ we have
\beq\label{AA0} \kk_{l,\bb}(t):=\inf_{|x-y|>l} \ff{K_1(t)V(x)+K_1(t) V(y)-2K_0(t)}{\bb^{-1}+ V(x)+V(y)}\in\R,\end{equation}
and when $K_1(t)>0$ and $l>0$ is large enough, $\kk_{l,\bb}(t)>0$.
Moreover,  $(H_4)$ and $(H_7)$ imply
\beq\label{AA} \beg{split} \aa_{l,\bb}(t):=& C_{\psi}\sup_{|x-y|\in (0,l)} \bigg\{\alpha_t \ff{|\nn V(x)-\nn V(y)| } {|x-y|\{\bb^{-1} + V(x)+V(y)\}}  \\
 & \quad + \ff{|\{\hat \si_t(x)-\hat\si_t(y)\}[(\hat \si_t(\cdot)^*\nn V)(x)+(\hat\si_t(\cdot)^*\nn V)(y)]|}{|x-y|\{\bb^{-1} + V(x)+V(y)\}}\bigg\}<\infty. \end{split}\end{equation}
  For   $K_0$, $\kk_{l,\bb}, \aa_{l,\bb}$ and $u_l$  given in $(H_7)$, $(H_8)$,   \eqref{AA0} and   \eqref{AA}, let
\beq\label{AA2} \ll_{l,\bb}(t):= \min\big\{\kk_{l,\bb}(t), \ u_l(t)- 2K_0(t)\bb -\aa_{l,\bb}(t)\big\},\ \ t\ge 0.\end{equation}

\beg{thm}\label{T8}   Assume  $(H_4)$, $(H_7)$ and $(H_8),$ with $\psi''\le 0$ when $\hat \si_t(\cdot)$ is non-constant.
 Then  $\eqref{E01}$ is well-posed for distributions in $\scr P_V$, and $P_t^*$
 satisfies
  \beq\label{EXP1}  \W_{\psi, \bb V}(P_t^*\mu, P_t^*\nu)\le \e^{-\int_0^t \{\ll_{l,\bb}(s) -\theta_s\}\d s }  \W_{\psi, \bb V}(\mu,\nu),\ \ t\ge 0, \mu,\nu\in \scr P_V.\end{equation}
  Consequently, if $(\si_t,b_t)$ is $t_0$-periodic   and
  $$\ll:= \int_0^{t_0} \{\ll_{l,\bb}(s) -\theta_s\}\d s >0,\ \ \int_0^{t_0} K_1(t)\d t>0,$$ then
    $\eqref{E01}$ has a unique invariant probability measure
$\bar\mu_0\in \scr P_V$ such that
\beq\label{EXP2}  \W_{\psi,\bb V}(P_{nt_0}^*\mu, \bar\mu_0)\le  \e^{-\ll n}  \W_{\psi, \bb V}(\mu,\bar\mu_0),\ \ n\in\mathbb N, \mu \in \scr P_V.\end{equation}
  \end{thm}

\beg{proof} The well-posedness and \eqref{EXP1} is included in   \cite[Theorem 2.1]{W21a}. So, it suffices to prove the existence of invariant probability measure $\bar\mu_0\in \scr P_V$ for the $t_0$-periodic case with $\ll>0$. Let $x_0\in \R^d$. By \eqref{EXP1} we have
$$\W_{\psi,\bb V} (P_{nt_0}^*\dd_{x_0}, P_{(n+m)t_0}^*\dd_{x_0})\le \e^{-\ll n} \W_{\psi,\bb V}(\dd_{x_0}, P_{m t_0}^*\dd_{x_0}),\ \ n,m\ge 1.$$
Therefore, it suffices to prove
\beq\label{PRW} \sup_{m\ge 1} \E V(X_{m t_0} )<\infty\ \text{for}\ X_0=x_0,\end{equation}
which together with the above inequality implies that $\{P_{n t_0}^*\dd_{x_0}\}_{n\ge 1}$ is a $\W_{\psi,\bb V}$-Cauchy sequence and its limit is
an invariant probability measure in $\scr P_V$. By \eqref{H120}, It\^o's formula and
$$\int_{m}^{(n+m)t_0} K_1(s)\d s= n \int_0^{t_0}K_1(s)\d s=:n \ll_0>0,\ \ n,m\in \mathbb N,$$  we obtain
we obtain
\beg{align*}\E V(X_{nt_0})&\le V(x_0)\e^{-\int_0^{nt_0} K_1(s)\d s} +\int_0^{nt_0} |K_0(s)| \e^{-\int_{s}^{nt_0}K_1(r)\d r}\d s\\
&\le V(x_0) +\sum_{i=0}^{n-1} \int_{i t_0}^{(i+1)t_0} C| K_0(s)| \e^{- \int_{(i+1)t_0}^{nt_0} K_1(r)\d r}\d s\\
&= V(x_0)+ \Big(\sum_{i=0}^{n-1} \e^{-(n-i-1)\ll_0}\Big)\int_0^{t_0}C |K_0(s)|\d s,\ \ n\ge 1,\end{align*}
which is bounded in $n\ge 1$  since $\ll_0:=\int_0^{t_0} K_1(t)\d t>0$. So, \eqref{PRW} holds.
\end{proof}
In the following example the SDE includes a class of  fully non-dissipative models, for instance when $\nn^{(1)}W\ge 0$,  in the sense that
$$\sup_{|x-y|=r} \<b_t(x,\mu)-b_t(y,\mu), x-y\> \ge 0,\ \ r>0, \mu\in \scr P.$$

\paragraph{ Example of Theorem 4.1.}

Let $\aa\in C([0,t_0];(0,\infty))$, $b_0\in C^1(\R^d)$ with $b_0(x)=- |x|^{p-1}x$ for $|x|\ge 1$, and $W_t\in C^2(\R^d\times\R^d)$ measurable in $t\in [0,t_0]$ with
\beq\label{KLN} \|\nn^{(1)}\nn ^{(2)}W_t\|_\infty+\|\nn^{(1)}W_t\|_\infty\le \vv \aa_t~~\mbox{and}~~\|\nn^{(1)}\nn ^{(1)}W_t\|_\infty\le \theta\alpha_t\end{equation} for some constant $\vv>0$.
 We take $t_0$-periodic $(b_t,\si_t)$ with
\beg{align*} &b_t(x,\mu) := \aa_t b_0(x) +\ff {\mu(\nn^{(1)}W_t(x,\cdot))} {1+\mu(V)},\\
&\si_t:=\ss{\aa_t} I_d,\ \ (t,x,\mu)\in [0,t_0]\times\R^d\times \scr P_V,\end{align*}
where  $V(x):=  \e^{|x|^{p}}$ for some $p\in [\ff 1 2,1].$  Moreover, let
\begin{equation}
 \tt \W_{V}(\mu, \nu):= \inf_{\pi\in \C(\mu, \nu)}\int_{\R^d\times \R^d}(1\wedge |x-y|)(1+V(x)+V(y))\pi(\d x, \d y)
\end{equation}
and
$$ \tt \W_{V}(P_t^*\mu, \bar\mu)\le c \e^{-\dd t}  \tt \W_{V}(\mu,\bar\mu),\ \ t\ge 0, ~\mu, \nu \in \scr P_V.$$
Then when $\vv>0$ is small enough,
there exist constants  $c,\ll>0$ such that
$$ \tt \W_{V}(P_{nt_0}^*\mu, \bar\mu)\le c \e^{-\ll n}  \tt \W_{V}(\mu,\bar\mu),\ \ n\in \mathbb N, \mu \in \scr P_V.$$

\beg{proof}   It is easy to see  that $(H_7)$ with
\beq\label{HH0} K_0=\aa_t \theta_0,\ \ K_1(t)= \aa_t \theta_1,\end{equation}    holds for some constants $\theta_0,\theta_1>0$,   $(H_8)$ holds for $\hat \si=0$. Next, let $D_0:=\|\nn b_0\|_\infty+\theta$ and let $l>0$ such that in \eqref{AA0}
\beq\label{HH0'}k_{l,\bb}(t):= \inf_{|x-y|\ge l} \frac{\theta_1 V(x)+\theta_1 V(y)- 2\theta_0}{\beta^{-1}+V(x)+V(y)}\ge k_0 \aa_t,\ \ t\in [0,t_0]\end{equation}
holds for some constant $k_0>0.$
Now, we take $\psi\in \Psi_l$ such that
$$2\psi''(r)+D_0 \psi'(r)\le -D_1\psi(r),\  \ r\in [0,l]$$ holds for some constant $D_1>0$,   for instance $\psi$ and $D_1$ ate  the first mixed eigenfunction and eigenvalue of $2\ff{\d^2}{\d r^2} +D_0 \ff{\d}{\d r}$ on $[0,l]$ with Dirichlet condition at $0$ and Neumann condition at $l$.

Then the first inequality in $(H_8)$ holds for
\beq\label{HH1} K_t:= \aa_t D_0, \ \ u_l(t):= D_1\aa_t,\ \ t\in [0,t_0].\end{equation}
Moreover, noting that $|V(x)-V(y)|\le c_0 \psi(|x-y|) (1+ V(x)+V(y))$ holds for some constant $c_0>0$,
we  find a constant $c_1>0$ such that
\beg{align*} |b(x,\mu)- b(x,\nu)|&\le \vv\Big(\ff{|\mu(\nn^{(1)}W(x,\cdot))- \nu(\nn^{(1)}W(x,\cdot))|}{1+\mu(V)\lor\nu(V)}+ \ff{\|\nn^{(1)}W\|_\infty |\mu(V)-\nu(V)|}{(1+\mu(V))(1+\nu(V))}\Big)\\
&\le \ff{ c_1 \vv \aa_t \{\W_{\psi,V}(\mu,\nu)  }{1+\mu(V)+\nu(V)}\le c_1 \vv \bb^{-1} \aa_t \hat{\W}_{\psi,\bb V}(\mu,\nu).\end{align*}
Combining this with $D_0:=\|\nn b^0\|_\infty +\theta$ and \eqref{KLN}, we obtain the second inequality in $(H_8)$ for
the above $K_t:= \aa_t D_0$ and
\beq\label{HH2}  \theta_t:= c_1\vv \bb^{-1}\aa_t,\ \ t\in [0,t_0].\end{equation}
Since  \eqref{AA} implies $\aa_{l,\bb}(t)\to 0$ as $\bb\to 0$, by \eqref{AA2}, \eqref{HH0}, \eqref{HH0'}, \eqref{HH1} and \eqref{HH2}, there exist  constants $\bb,\vv_0>0,k_1$ such that for any $ \vv\in (0,\vv_0]$
$$\ll_{l,\bb}(t)- \theta_t \ge k_1,\ \ t\in [0,t_0].$$
 Then the desired assertion follows from Theorem \ref{T8} and the fact that
 $$ C^{-1} \tt \W_V \le \W_{\psi,\bb V} \le C\tt \W_V$$
holds for some constant $C>1$. \end{proof}

\section{Extensions to reflecting McKean-Vlasov SDEs }

In this section,  we investigate the exponential ergodicity for the reflecting McKean-Vlasov SDE \eqref{E1} on a convex domain $D$.
By the convexity, the reflection on boundary does not make any trouble in the proofs of previous results on ergodicity, so that all these results
work also for \eqref{E1}.

Let $T\pp D$ be the tangent space of $\pp D$, which is well defined when $\pp D$ is $C^1$.

\beg{thm}\label{T6} Let $D$ be convex, $b,\si\in C([0,\infty)\times\bar D\times \scr P_2(\bar D))$, and  in $(H_1)$ -$(H_3)$ we use  $(\bar D,\scr P_2(\bar D))$ to replace $(\R^d,\scr P_2)$, and in $(H_3)(2)$ assume further that $\pp D$ is $C^2$ and there exists a measurable function
$h: [0,\infty)\times\pp D\to [0,\infty)$ such that
\beq\label{**X} \<\{\nn_{n}(\si_t\si_t^*)\} v,v\>|_{\pp D}\ge 0,\ \ (\si_t\si_t^*  v - h_t v)|_{\pp D}=0,\ \ v\in T\pp D, t\ge 0.\end{equation}
Then assertions in Theorem $\ref{T1.1}$ holds for $\eqref{E1}$ replacing $\eqref{E01}$. \end{thm}

\beg{proof} By \cite[Theorem 2.6]{W21b}, $(H_1)$ implies that \eqref{E1} is well-posed for distributions in $\scr P_2(\bar D)$ and satisfies
\beq\label{SOP2} \W_2(P_t^*\mu,P_t^*\nu)^2\le \e^{\int_0^t (K_1(s)+K_2(s))\d s} \W_2(\mu,\nu),\ \ \mu,\nu\in \scr P_2(\bar D).\end{equation}
Let   $x_0\in D.$ Since $D$ is convex, we have
$\<x-x_0, \n(x)\>\le 0$ for $x\in \pp D$, so that as in the proof of Theorem \ref{T1.1}, by $(H_1)$ and applying It\^o's formula  to $|X_t-x_0|^2$ for $X_0=x_0,$ we obtain
$$\d |X_t-x_0|^2\le \Big\{c+ \Big(K_1(t)+\ff{\ll}{2t_0}\Big) |X_t-x_0|^2+K_2(t) E|X_t-x_0|^2\Big\}\d t +\d M_t$$
for some martingale $M_t$. Since $\ll>0$, this and the proof leading to \eqref{SOP} gives the same estimate, so that by \eqref{SOP2} we prove the first assertion.

Under $(H_2)$ holds, by \cite[Theorem 2.4]{W21b}, there exists a constant $c_1>0$ such that
\beq\label{SOP3} \W_2(P_{t_0}^*\mu, P_{t_0}^*\nu)\le c_1 \H(\mu|\nu),\ \ \mu,\nu\in \scr P_2(\bar D).\end{equation}
So, as in the proof of Theorem \ref{T1.1}, it remains to  prove the  Talagrand inequality
\beq\label{TTI'} \W_2(\mu,\bar\mu_0)^2\le c_2 \Ent(\mu|\bar\mu_0),\ \ \mu\in \scr P_2(\bar D)\end{equation}  for some constant $c_2>0$.

When $(H_3)$(1) holds, by the convexity of $D$, for $(\bar X_t^x, \bar X_t^y)$ solving the following SDE  with $\bar X_0^x=x, \bar X_0^y\in \bar D$:
 $$\d \bar X_t= b_t(\bar X_t, \bar \mu_0)\d t+\si_t\d W_t+\n(\bar X_t)\d l_t,$$
 $(H_1)$ implies
$$\d |\bar X_t^x-\bar X_t^y|^2\le  K_1(t)|\bar X_t^x-\bar X_t^y|^2 \d t,\ \ t\ge 0,$$
so that
$$ |\bar X_t^x-\bar X_t^y|^2\le \e^{\int_0^t K_1(s) \d s}|x-y|^2,\  \ x,y\in \bar D,\ t\ge 0.$$
Thus, the associated $\bar P_t$ satisfies the gradient estimate \eqref{GRD}. Then \eqref{TTI'} holds as shown in the proof of Theorem \ref{T1.1},

When $(H_3)(2)$ holds, the corresponding proof in that of Theorem \ref{T1.1} also works provided Lemma \ref{LN1}(2) holds for \eqref{E1}.
According to its proof it suffices to prove \cite[Theorem 4.1]{CM} for \eqref{E1}, which is included in the following Lemma \ref{LN2}.
\end{proof}
Let $\GG_t^1$  and $\GG_2^t$ be in $\eqref{GGG}$ for $\si_t,b_t$ not depending on $\mu$ on a convex $C^2$ domain $\bar D$ replacing $\R^d$, where
$b_t(x)$ is $C^1$ in $x$,   $\si_t(x)$ is $C^1$ in $t$ and $C^2$ in $x$.
Consider
the reflecting SDE
\beq\label{E1'} \d X_{s,t}= b_t(X_{s,t})\d t + \si_t(X_{s,t})\d W_t+\n(X_t)\d l_t,\ \ t\ge s.\end{equation}
Let $P_{s,t}^*\mu=\L_{X_{s,t}}$ for the solution with $\L_{X_{s,s}}=\mu.$ The generator is
$$L_t := \ff 1 2 {\rm tr}\{\si_t\si_t^*\nn^2\} + b_t\cdot\nn,\ \ t\ge 0.$$
We have the following lemma, which extends \cite[Theorem 4.1]{CM} to the reflecting case.

\beg{lem}\label{LN2} Let $\{\GG_t^i\}_{i=1,2, t\ge 0}$ be in $\eqref{GGG}$ on a convex $C^2$ domain $\bar D$ replacing $\R^d$  for $\si_t,b_t$ not depending on $\mu$, and let \eqref{**X} hold.
Let $\gg\in L_{loc}^1([0,\infty);\R)$   such that
\beq\label{GGO} \GG_t^2(f,f)   \ge \gg_t \GG_t^1(f,f),\ \ f\in C^3(\bar D).\end{equation}
Let $s\ge 0,q_s>0$ and $\nu_s\in \scr P(\bar D).$ If the log-Sobolev inequality
\beq\label{VS}\nu_s(f^2\log f^2)\le 4q_s \nu_s(\GG_1^s(f,f)),\ \ f\in C^1_b(\bar D), \nu_s(f^2)=1\end{equation} holds,  then  for any $t>s$,
$\nu_t:=P_{s,t}^*\nu_s$   satisfies
\beq\label{VT}\nu_t(f^2\log f^2)\le 4q_t \nu_t(\GG_1^t(f,f)),\ \ f\in C^1_b(\bar D), \nu_t(f^2)=1,\end{equation}  where
\beq\label{QT}q_t:= q_s\e^{-2\int_s^t \gg_r\d r}+\int_s^t\e^{-2\int_r^t\gg_u\d u}\d r, \ \ t\ge s.\end{equation}
\end{lem}

\beg{proof} Let $P_{s,t}f(x):= (P_{s,t}^*\dd_x)(f).$ We first  prove
\beq\label{GRD1} \ss{\GG_1^{s_0}(P_{s_0,s_1}f,P_{s_0,s_1}f) }\le \e^{-\int_{s_0}^{s_1} \gg_t\d t} P_{s_0,s_1} \ss{\GG_1^{s_1}(f,f)},\ \ s_1\ge s_0\ge 0.\end{equation}
for $f\in C_b^\infty(\bar D).$

 By the Bochner-Weitzenb\"ock formula,     the inequality \eqref{GGO} is equivalent to
\beq\label{GGO'} \GG_t^2(f,f)   \ge \gg_t \GG_t^1(f,f) + \ff{|\nn \GG_1^t(f,f)|^2}{4 \GG_1^t(f,f)},\ \ f\in C^3(\R^d).\end{equation}
Next, since $\pp D$ is $C^2$ and convex, the second fundamental form is non-negative, i.e.
$$\II(\nn f,\nn f)(x):= -\<\nn_{\nn f}\n, \nn f\>(x)=\Hess_f(\n, \nn f)(x)\ge 0,\ \ f\in C^2(\bar D), Nf|_{\pp D}=0, x\in \pp D,$$
where the second equality follows from $\nn_{\nn f} \<\n, \nn f\>|_{\pp D}=0$ due to $\<\n,\nn f\>|_{\pp D}=0.$
Combining this with \eqref{**X},  we see that  for any $f\in C_b^\infty(\bar D), s_1\ge t\ge 0,$ and $x\in \pp D,$
\beg{align*} &\big\<\n, \nn \GG_1^s(P_{t,s_1} f^2, P_{t,s_1} f^2)\big\>(x)\\
&= \big\<\{\nn_\n (\si_t\si_t^*)\}\nn P_{t,s_1} f^2, \nn P_{t,s_1} f^2\big\>(x) + 2 \Hess_{P_{t,s_1}f^2} (\n, (\si_t\si_t^*)\nn P_{t,s_1} f^2\big\>(x)\ge 0.\end{align*}
 Combining this with \eqref{GGO'} and  applying It\^o's formula to \eqref{E1'}, for any $s_1>s_0\ge 0$,
 \beg{align*} &\d \ss{\GG_1^t (P_{t, s_1} f^2, P_{t,s_1} f^2)}(X_{s_0,t})\\
 &= \bigg\{\ff{\ff 1 2 (\pp_t \GG_1^t)(P_{t, s_1} f^2, P_{t,s_1} f^2)- \GG_1^t (L_tP_{t, s_1} f^2, P_{t,s_1} f^2)}{\ss{\GG_1^t (P_{t, s_1} f^2, P_{t,s_1} f^2)}} + L_t \ss{\GG_1^t (P_{t, s_1} f^2, P_{t,s_1} f^2)} \bigg\}(X_{s_0,t})\d t \\
 &\quad +\d M_t+ \Big\<\n, \nn \ss{\GG_1^t(P_{t,s_1} f^2, P_{t,s_1} f^2)\big\>}(X_{s_0,t})\Big\>\d l_t\\
 &\ge \bigg\{\ff{\GG_t^2(P_{t, s_1} f^2, P_{t,s_1} f^2)}{\GG_t^1(P_{t, s_1} f^2, P_{t,s_1} f^2)^{\ff 1 2}}
 - \ff{|\nn \GG_1^sP_{t, s_1} f^2, P_{t,s_1} f^2)|^2}{4\GG_1^t(P_{t, s_1} f^2, P_{t,s_1} f^2)^{\ff 3 2}} \bigg\}(X_{s_0,t})\d t+\d M_t\\
 &\ge \gg_t \ss{\GG_1^t (P_{t, s_1} f^2, P_{t,s_1} f^2)}(X_{s_0,t})\d t +\d M_t,\ \ t\in [s_0,s_1]\end{align*}
holds for some martingale $(M_t)_{t\in [s_0,s_1]}$. By Gronwall's lemma this implies \eqref{GRD1}.

By \eqref{GRD1}, the desired assertion follows from a standard semigroup argument, we include below for completeness. Let $f\in C_b^2(\bar D)$ with $\inf f^2>0$.
By the chain rule and Schwarz inequality, \eqref{GRD1} implies
\beq\label{GRD2}\beg{split} &\GG_1^s\big(\ss{P_{s,t}f^2},\ss{P_{s,t}f^2}\big)= \ff {\GG_1^s(P_{s,t}f^2, P_{s,t}f^2 )}{4P_{s,t}f^2}\\
&\le \ff{\e^{-2\int_s^t \gg_r\d r}(P_{s,t}\ss{\GG_1^t(f^2,f^2)})^2}{4 P_{s,t}f^2} \le \e^{-2\int_s^t\gg_r\d r} P_{s,t}\GG_1^t(f,f),\ \ t\ge s\ge 0.\end{split}\end{equation}
So,
\beg{align*} & P_{s,t} (f^2\log f^2) - (P_{s,t}f^2)\log P_{s,t} f^2= \int_s^t \ff{\d}{\d r} P_{s,r} \big\{(P_{r,t}f^2)\log (P_{r,t}f^2)\big\}\d r \\
&=\int_s^t P_{s,r}\ff{\GG_1^r(P_{r,t}f^2, P_{r,t}f^2)}{P_{r,t}f^2}\d r\le 4 (P_{s,t} \GG_1^t(f,f)) \int_s^t \e^{-2\int_r^t \gg_u\d u}\d r.\end{align*}
Combining this with \eqref{VS} and \eqref{GRD2}, we obtain
\beg{align*} &\nu_t(f^2\log f^2)= \nu_s(P_{s,t}(f^2\log f^2))\\
&\le 4 \nu_s\big(P_{s,t}\GG_1^t(f,f)) \int_s^t \e^{-2\int_r^t\gg_u\d u}\d r+\nu_s\big((P_{s,t}f^2)\log (P_{s,t}f^2)\big) \\
&\le 4\nu_t\big(\GG_1^t(f,f)\big) \int_s^t \e^{-2\int_r^t\gg_u\d u}\d r + 4 q_s\nu_s\Big(\GG_1^s\Big(\ss{P_{s,t}f^2},\ss{P_{s,t}f^2}\Big)\Big)
+\big(\nu_s(P_{s,t}f^2)\big)\log\big(\nu_s(P_{s,t}f^2)\big) \\
&\le 4 \nu_t\big(\GG_1^t(f,f)\big) \bigg( \int_s^t \e^{-2\int_r^t\gg_u\d u}\d r+ q_s \e^{-2\int_s^t \gg_u\d u}\bigg\}+ \nu_t(f^2)\log \nu_t(f^2).\end{align*}
Therefore, \eqref{VT} holds for $q_t$ in \eqref{QT}.
\end{proof}

 Finally, we have the following extensions of Theorems \ref{T3-1} and \eqref{T8}.

 \beg{thm} Let $D$ be convex, use $\bar D$ replace $\R^d$ in $(H_4)$-$(H_8)$,  and  in $(H_7)$ we assume further  $\<\nn V,\n\>|_{\pp D}\le 0$. Then  assertions in    Theorem $\ref{T3-1}$   and Theorem $\ref{T8}$ hold for $\eqref{E1}$  replacing $\eqref{E01}$.
 \end{thm}

 \beg{proof} The well-posedness of \eqref{E1} as well as estimates \eqref{EXP1'0} and \eqref{EXP1}
 have been included in \cite[Theorems 2.7, 2.8]{W21b}, so that the other assertions follow from the proofs of Theorems \ref{T3-1} and \ref{T8}.
  Indeed, the proof of Theorem \ref{T3-1} has nothing to do with the reflection. Moreover,    by It\^o's formula, \eqref{H120} and    $\<\n,\nn V\>|_{\pp D}\le 0$,   we derive
 $$\d V(X_t)\le \{K_0(t)-K_1(t)V(X_t)\}\d t+\d M_t$$ for some local martingale $M_t$, so that the proof of \eqref{PRW} works also for the present case.

 \end{proof}


\end{document}